# ON OVERLOAD IN A STORAGE MODEL, WITH A SELF-SIMILAR AND INFINITELY DIVISIBLE INPUT

By J. M. P. Albin[1] and Gennady Samorodnitsky[2]

*Chalmers University of Technology and Cornell University*

Let $\{X(t)\}_{t \geq 0}$ be a locally bounded and infinitely divisible stochastic process, with no Gaussian component, that is self-similar with index $H > 0$. Pick constants $\gamma > H$ and $c > 0$. Let $\nu$ be the Lévy measure on $\mathbb{R}^{[0,\infty)}$ of $X$, and suppose that $R(u) \equiv \nu(\{y \in \mathbb{R}^{[0,\infty)} : \sup_{t \geq 0} y(t)/(1 + ct^\gamma) > u\})$ is suitably "heavy tailed" as $u \to \infty$ (e.g., subexponential with positive decrease). For the "storage process" $Y(t) \equiv \sup_{s \geq t}(X(s) - X(t) - c(s - t)^\gamma)$, we show that $\mathbf{P}\{\sup_{s \in [0, t(u)]} Y(s) > u\} \sim \mathbf{P}\{Y(\hat{t}(u)) > u\}$ as $u \to \infty$, when $0 \leq \hat{t}(u) \leq t(u)$ do not grow too fast with $u$ [e.g., $t(u) = o(u^{1/\gamma})$].

**1. Introduction.** Let $X = \{X(t)\}_{t \geq 0}$ be an infinitely divisible (i.d.) stochastic process, with no Gaussian component, that is self-similar with index $H > 0$
($H$-s.s.).

Given constants $c > 0$ and $\gamma > H$, we consider the *storage process*

$$(1.1) \qquad Y(t) = \sup_{s \geq t}(X(s) - X(t) - c(s - t)^\gamma) \qquad \text{for } t \geq 0.$$

Intuitively, an $H$-s.s. process grows as $t^H$ with time $t$, and so $\gamma > H$ should make $Y$ finite valued. Nevertheless, this is not so in general (see Example 2). The assumptions in our theorems will, however, ensure such finiteness of $Y$. The reason for the name "storage process" comes from the case $\gamma = 1$, with $X(t)$ denoting the total inflow into a storage facility by time $t$, and $c$ the (demand) rate at which stock at the facility is depleted; then $Y(t)$ tells how much extra storage capacity one will need in the future over what is

Received July 2002; revised May 2003.
[1]Supported by NFR Grant M 650-19981841/2000, and by STINT Grant Dnr 00/28.
[2]Supported by NSF Grant DMS-00-71073 at Cornell University.
*AMS 2000 subject classifications.* Primary 60G18, 60G70; secondary 60E07, 60G10.
*Key words and phrases.* Heavy tails, infinitely divisible process, Lévy process, self-similar process, stable process, stationary increment process, subexponential distribution, storage process.

being used at time $t$. For an input process $X$ with stationary increments (s.i.), the storage process $Y$ is stationary (if finite).

The process $Y$ has been used in financial applications under the name of "drawdown" [e.g., Dacorogna, Gençay, Müller and Pictet (2001)], and is important in queueing applications; for example, to model teletraffic, when $X$ is Gaussian $H$-s.s.s.i. (i.e., fractional Brownian motion), with $H \geq \frac{1}{2}$ and linear service $\gamma = 1$ [e.g., Norros (1994) and Piterbarg (2001)]. In this case, building on Hüsler and Piterbarg (1999), Piterbarg [(2001), Theorem 5] gave a version of the remarkable property (1.2), that was a triggering influence for us. Recognizing this, we name that property after him. However, the Gaussian problem Piterbarg studied is very different from ours with i.d. processes, and his proof, by Gaussian field theory, does not relate to non-Gaussian settings.

In Section 3 we study the probability for overload during a time interval,

$$\mathbf{P}\left\{ \sup_{s \in [0,t]} Y(s) > u \right\} \qquad \text{as } u \to \infty.$$

If for each choice of a constant $t > 0$, it holds that

(1.2a)    $$\lim_{u \to \infty} \frac{\mathbf{P}\{\sup_{s \in [0,t]} Y(s) > u\}}{\mathbf{P}\{Y(\hat{t}(u)) > u\}} = 1 \qquad \text{whenever } 0 \leq \hat{t}(u) \leq t,$$

then we say that the process $Y$ has the *Piterbarg property*. The similar statement, for which the length $t = t(u)$ of the interval may depend on the level $u$,

(1.2b)    $$\lim_{u \to \infty} \frac{\mathbf{P}\{\sup_{s \in [0,t(u)]} Y(s) > u\}}{\mathbf{P}\{Y(\hat{t}(u)) > u\}} = 1 \qquad \text{whenever } 0 \leq \hat{t}(u) \leq t(u),$$

will be referred to as the *generalized Piterbarg property*.

One indication of the unusual behavior of $Y$ is that (1.2) implies, as $u \to \infty$,

$$\mathbf{P}\left\{ \bigcap_{i=1}^{n} \{Y(\hat{t}_i) > u\} \,\Big|\, \sup_{s \in [0,t]} Y(s) > u \right\}$$

$$\geq 1 - \sum_{i=1}^{n} \left( 1 - \frac{\mathbf{P}\{Y(\hat{t}_i) > u\}}{\mathbf{P}\{\sup_{s \in [0,t]} Y(s) > u\}} \right) \to 1.$$

Thus overload periods within $[0,t]$ are long enough to include any $\hat{t}_1, \ldots, \hat{t}_n \in [0,t]$.

This last conclusion leads us naturally to the question whether one can replace the minimum $\bigwedge_{i=1}^{n} Y(\hat{t}_i)$ taken over a finite collection of points in

$[0, t]$ by the infimum over the entire interval. That is, we would like to know if

$$(1.2c) \qquad \lim_{u \to \infty} \frac{\mathbf{P}\{\sup_{s \in [0,t]} Y(s) > u\}}{\mathbf{P}\{\inf_{s \in [0,t]} Y(s) > u\}} = 1.$$

This we call the *strong Piterbarg property*, whether or not $t$ is a function of $u$.

With $\nu$ being the Lévy measure on $\mathbb{R}^{[0,\infty)}$ of $X$ (see Section 2.3), denote

$$(1.3) \quad R(u) \equiv 1 \wedge \nu\left(\left\{y \in \mathbb{R}^{(0,\infty) \cap \mathbb{Q}} : \sup_{t \in (0,\infty) \cap \mathbb{Q}} \frac{y(t)}{1 + ct^\gamma} > u\right\}\right) \qquad \text{for } u \in \mathbb{R}.$$

We will make assumptions about "heavy tails" for the function $R$ (e.g., subexponentiality together with positive decrease; see Section 2.1). Under additional technical assumptions on $X$, we establish the generalized Piterbarg property, when $t(u)$ does not grow too fast with $u$ [e.g., $t(u) = o(u^{1/\gamma})$]. Under the same assumptions we will show that the strong Piterbarg property holds as well. Under certain weaker assumptions, we prove a weaker so-called $O$-version of (1.2), that is, that the probability ratios in (1.2) are bounded away from zero and infinity.

Our main "external tool" in proofs is Theorem 2.1 on subexponential functionals of i.d. processes by Rosiński and Samorodnitsky (1993); see Section 2.3.

The Piterbarg properties are quite unusual. For example, only a degenerate $\alpha$-stable or Gaussian process $Y$ can have them; see Example 7.

In Section 4, we give a discussion, with examples of application, and counterexamples, for i.d. $H$-s.s. processes $X$ given as stochastic integrals with respect to heavy-tailed i.d. random measures (see Section 2.3). This includes $\alpha$-stable processes.

**2. Classes of functions and stochastic processes.** It will be convenient to devote a separate section to describe classes of functions and stochastic processes, that feature in the rest of the article. In addition, some basic relations between these classes, and some important representation properties, are listed for easy reference.

2.1. *Classes of functions.* In this section, $f : \mathbb{R} \to (0, \infty)$ denotes a nonincreasing function with $\lim_{u \to \infty} f(u) = 0$, and $g : \mathbb{R} \to (0, \infty)$ a measurable function.

The function $f$ is *subexponential*, $f \in \mathcal{S}$, if there exist independent identically distributed random variables $\xi$ and $\eta$, such that

$$f(u) \sim \mathbf{P}\{\xi > u\} \quad \text{and} \quad \mathbf{P}\{\xi + \eta > u\} \sim 2\mathbf{P}\{\xi > u\} \qquad \text{as } u \to \infty.$$

The function $f$ is *O-regularly varying*, $f \in OR$, if

$$\liminf_{u \to \infty} \frac{f(\lambda u)}{f(u)} > 0 \qquad \text{for some } \lambda > 1.$$

The function $f$ has *positive decrease*, $f \in PD$, if

$$\limsup_{u \to \infty} \frac{f(\lambda u)}{f(u)} < 1 \qquad \text{for some } \lambda > 1.$$

The function $f$ is *extended regularly varying*, $f \in ER$, if

$$\liminf_{u \to \infty} \frac{f(\lambda u)}{f(u)} \geq \lambda^{-b} \qquad \text{for } \lambda \geq 1, \text{for some constant } b \geq 0.$$

Note.  The definitions of $OR$, $PD$ and $ER$ are more complicated than those given above for a general nonmonotone $f$.

The function $g$ is *regularly varying with index* $\rho \in \mathbb{R}$, $g \in RV(\rho)$, if

$$\lim_{u \to \infty} \frac{g(\lambda u)}{g(u)} = \lambda^{\rho} \qquad \text{for } \lambda \geq 1 \text{ (or, equivalently, for } \lambda > 0).$$

Here the convergence must, in fact, be locally uniform.

Notice that the function $g \circ \log$ belongs to $RV(0)$, which we denote $g \in \mathcal{L}$, if

$$\lim_{u \to \infty} \frac{g(u + \lambda)}{g(u)} = 1 \qquad \text{for } \lambda \geq 0 \text{ (or, equivalently, for } \lambda \in \mathbb{R}).$$

We have $ER \cap PD \subseteq OR \cap PD \cap \mathcal{S}$, and a monotone $f \in \bigcup_{\rho < 0} RV(\rho)$ belongs to all these classes. Further, $OR \cap \mathcal{L} \subseteq \mathcal{S} \subseteq \mathcal{L}$.

The classes of functions above, and the listed relations between them, are well known from the literature. See, for example, Bingham, Goldie and Teugels (1987).

2.2. *Classes of stochastic processes.*  In the remainder of this article, $X = \{X(t)\}_{t \geq 0}$ denotes a separable stochastic process, that is continuous in probability and locally bounded (bounded on any given compact interval) a.s., and is defined on a complete probability space $(\Omega, \mathfrak{F}, \mathbf{P})$. We refer to these requirements as *Condition* $\mathfrak{X}$. Depending on the context, further requirements on $X$ will be imposed later.

We write $\{\tilde{X}(t)\}_{t \geq 0} \overset{d}{=} X$ when the finite-dimensional distributions (f.d.d.'s) of the processes $\tilde{X}$ and $X$ agree. For example, $X$ is stationary, if $X(\cdot + h) \overset{d}{=} X$ for $h \geq 0$.

The process $X$ has *stationary increments*, if

$$X(\cdot + h) - X(h) \overset{d}{=} X - X(0) \qquad \text{for each } h \geq 0.$$

The process $X$ is *self-similar with index* $H > 0$ ($H$-s.s.), if

$$a^{-H} X(a\cdot) \stackrel{d}{=} X \qquad \text{for each } a > 0.$$

The process $X$ is *infinitely divisible* (i.d.), if for each $n \in \mathbb{N}$, there exist independent processes $\{\tilde{X}_1(t)\}_{t\geq 0}, \ldots, \{\tilde{X}_n(t)\}_{t\geq 0}$, such that

$$\tilde{X}_1 \stackrel{d}{=} \cdots \stackrel{d}{=} \tilde{X}_n \quad \text{and} \quad \tilde{X}_1 + \cdots + \tilde{X}_n \stackrel{d}{=} X.$$

The process $X$ is $\alpha$-*stable*, $\alpha \in (0, 2]$, if for each $n \in \mathbb{N}$, there exists a constant process $C_n$, such that, taking independent copies $\{\tilde{X}_k\}_{k=1}^n$ of $X$,

$$n^{-1/\alpha}(\tilde{X}_1 + \cdots + \tilde{X}_n) + C_n \stackrel{d}{=} X.$$

In particular, it turns out, the process $X$ is Gaussian if and only if it is two-stable.

The process $X$ is *strictly* $\alpha$-*stable* if, taking independent copies $\{\tilde{X}_k\}_{k=1}^\infty$ of $X$,

$$n^{-1/\alpha}(\tilde{X}_1 + \cdots + \tilde{X}_n) \stackrel{d}{=} X \qquad \text{for } n \in \mathbb{N}.$$

A process $X$ is $H$-s.s. if and only if the Lamperti transformed process $e^{-H\cdot}X(e^\cdot)$ is stationary [Lamperti (1962)].

$\alpha$-stable processes are i.d. Clearly, an $\alpha$-stable process $X$ is strictly $\alpha$-stable if it is symmetric $\alpha$-stable ($S\alpha S$) ($\alpha$-stable with $X \stackrel{d}{=} -X$).

Of course, the classes of processes mentioned above are all quite basic, as are the indicated relations between. See, for example, Samorodnitsky and Taqqu (1994) for further information, and for an extensive bibliography.

2.3. *I.d. stochastic processes.* The f.d.d.'s of an i.d. process $\{X(t)\}_{t\in T}$, $T = [0, \infty)$, with no Gaussian component, can be described by means of a *Lévy measure* $\nu$ on the cylindrical $\sigma$-algebra $\mathcal{B}$ on $\mathbb{R}^T$, and a *localization* parameter $\mu \in \mathbb{R}^T$.

Let $\pi_\tau$ be the projection of $\mathbb{R}^T$ on $\mathbb{R}^\tau$, and let $\mathcal{B}_\tau$ be the Borel sets in $\mathbb{R}^\tau$, for $\tau \in \mathcal{T} \equiv \{\tau \subseteq T : 1 \leq \#\tau < \infty\}$. According to Maruyama (1970), a measure $\nu$ on $\mathcal{B}$ is a Lévy measure for $X$, if $\nu \circ \pi_\tau^{-1}$ is a Lévy measure on $\mathcal{B}_\tau$ [i.e, if $1 \wedge |\cdot|^2 \in \mathbb{L}^1(\mathbb{R}^\tau, \nu \circ \pi_\tau^{-1})$] for each $\tau \in \mathcal{T}$, and there exists a $\mu \in \mathbb{R}^T$ such that

$$(2.1) \quad \mathbf{E}\{e^{i\langle\theta, X\rangle}\} = \exp\left\{i\langle\theta, \mu\rangle + \int_{\mathbb{R}^T} (e^{i\langle\theta, x\rangle} - 1 - i\langle\theta, \kappa(x)\rangle)\, d\nu(x)\right\}$$
$$\text{for } \theta \in \mathbb{R}^{(T)}.$$

Here we use the notation

$$\mathbb{R}^{(T)} = \{x \in \mathbb{R}^T : \#\{t \in T : x(t) \neq 0\} < \infty\},$$

$$\langle x, y\rangle = \sum_{t\in T} x(t)y(t) \qquad \text{for } x \in \mathbb{R}^{(T)} \text{ and } y \in \mathbb{R}^T,$$

$$\kappa(x)(t) = x(t)\mathbf{1}_{[-1,1]}(|x(t)|) \qquad \text{for } x \in \mathbb{R}^T \text{ and } t \in T.$$

[A general i.d. $X$ can be represented as $X \overset{d}{=} X_1 + X_2$, with $X_1$ and $X_2$ independent, $X_1$ i.d. with no Gaussian component as in (2.1), and $X_2$ zero-mean Gaussian

$$\mathbf{E}\{e^{i\langle\theta,X_2\rangle}\} = \exp\left\{-\tfrac{1}{2} \sum_{s,t\in T,\theta(s),\theta(t)\neq 0} \theta(s)\theta(t)\mathbf{E}\{X_2(s)X_2(t)\}\right\}$$
$$\text{for } \theta \in \mathbb{R}^{(T)}.]$$

We now turn to the task of constructing and representing i.d. processes.

Let $(S, \mathfrak{S}, \lambda)$ be a $\sigma$-finite measure space, and put $\mathfrak{S}_0 \equiv \{A \in \mathfrak{S} : \lambda(A) < \infty\}$. An (independently scattered) i.d. *random measure* (with no Gaussian component), with control measure $\lambda$, is a map $M : \mathfrak{S}_0 \to \mathbb{L}^0(\Omega, \mathfrak{F})$ such that, for $A \in \mathfrak{S}_0$,

$$(2.2) \qquad \mathbf{E}\{e^{i\theta M(A)}\} = \exp\left\{\int_A \left(i\theta m + \int_{\mathbb{R}} (e^{i\theta x} - 1 - i\theta\kappa(x))\, \rho(\cdot, dx)\right) d\lambda\right\}$$
$$\text{for } \theta \in \mathbb{R}.$$

Here the localization $m \in \mathbb{L}^0(S)$ satisfies $\{\mathbf{1}_A m\}_{A\in\mathfrak{S}_0} \subseteq \mathbb{L}^1(S, \lambda)$, while $\rho(s, \cdot)$ is a Lévy measure on $\mathbb{R}$ for $s \in S$, such that $\rho(\cdot, B) \in \mathbb{L}^0(S)$ for Borel sets $B \subseteq \mathbb{R}$, and

$$(2.3) \quad F(A \times \cdot) \equiv \int_A \rho(s, \cdot)\lambda(ds) \text{ is a Lévy measure on } \mathbb{R} \text{ for each } A \in \mathfrak{S}_0.$$

The stochastic integral $\int_S f\, dM$ is well defined (in a $\mathbf{P}$-sense), for $f \in \mathbb{L}^0(S)$ with

$$(2.4) \qquad \begin{aligned} &\int_S \int_{\mathbb{R}} (1 \wedge |xf|^2)\rho(\cdot, dx)\, d\lambda \\ &\quad \vee \int_S \left| mf + \int_{\mathbb{R}} (\kappa(xf) - f\kappa(x))\rho(\cdot, dx) \right| d\lambda < \infty \end{aligned}$$

[Rajput and Rosiński (1989), Section 2]. In that case, $\int_S f\, dM$ is i.d., with

$$(2.5) \qquad \begin{aligned} &\mathbf{E}\left\{\exp\left(i\theta \int_S f\, dM\right)\right\} \\ &= \exp\left\{\int_S \left(i\theta m f + \int_{\mathbb{R}} (e^{i\theta x f} - 1 - i\theta\kappa(x)f)\rho(\cdot, dx)\right) d\lambda\right\}. \end{aligned}$$

In the language of (2.1), the (process consisting of a) single i.d. random variable $\int_S f\, dM$ has Lévy measure $\nu$ on $\mathbb{R}$, and localization parameter $\mu \in \mathbb{R}$, given by

$$\nu(B) = F(\{(s, x) \in S \times \mathbb{R} : xf(s) \in B\}),$$
$$\mu = \int_S \left(mf + \int_{\mathbb{R}} (\kappa(xf) - \kappa(x)f)\rho(\cdot, dx)\right) d\lambda.$$

In particular, for example, by Feller [(1971), page 571], $\int_S f \, dM$ is nonnegative, if and only if

$$xf(s) \geq 0 \qquad \text{a.e. } (F),$$

$$(2.6) \qquad \int_S \left( mf - \int_{\mathbb{R}} \kappa(x) f \rho(\cdot, dx) \right) d\lambda \geq 0,$$

$$\int_S \int_{\mathbb{R}} (1 \wedge |xf|) \rho(\cdot, dx) \, d\lambda < \infty.$$

Pick $f_t \in \mathbb{L}^0(S)$ satisfying (2.4) for $t \geq 0$. The following process is i.d.:

$$(2.7) \qquad X \stackrel{d}{=} \left\{ \int_S f_t \, dM \right\}_{t \geq 0}$$

with Lévy measure in (2.1) given by $\nu = F \circ T_f^{-1}$, where $S \times \mathbb{R} \ni (s, x) \mapsto T_f(s, x) = x f_{(\cdot)}(s) \in \mathbb{R}^{[0, \infty)}$. With $\mathbb{Q}^+ = (0, \infty) \cap \mathbb{Q}$, the function $R$ in (1.3) thus satisfies

$$(2.8) \qquad R(u) = 1 \wedge \int_S \rho\left( s, \mathbb{R} \setminus \left[ \frac{-u}{\sup_{t \in \mathbb{Q}^+} f_t(s)^- / (1 + ct^\gamma)}, \right. \right.$$
$$\left. \left. \frac{u}{\sup_{t \in \mathbb{Q}^+} f_t(s)^+ / (1 + ct^\gamma)} \right] \right) d\lambda(s).$$

By (2.5), the process $X$ in (2.7) is $H$-s.s., if and only if

$$(2.9) \qquad \int_S \left( ima^{-H} \sum_{j=1}^n \theta_j f_{at_j} \right.$$
$$\left. + \int_{\mathbb{R}} \left( \exp\left( ixa^{-H} \sum_{j=1}^n \theta_j f_{at_j} \right) - 1 - i\kappa(x) a^{-H} \sum_{j=1}^n \theta f_{at_j} \right) \rho(\cdot, dx) \right) d\lambda$$

does not depend on $a > 0$ for any choice of $n \in \mathbb{N}$, $t_1, \ldots, t_n \geq 0$ and $\theta_1, \ldots, \theta_n \in \mathbb{R}$. Similarly, $X$ is $H$-s.s.s.i., if and only if $X(0) = 0$, and, with obvious notation,

$$\int_S \left( ima^{-H} \langle \theta, f_{at+h} - f_h \rangle \right.$$
$$\left. + \int_{\mathbb{R}} \left( e^{ixa^{-H} \langle \theta, f_{at+h} - f_h \rangle} - 1 - i\kappa(x) a^{-H} \langle \theta, f_{at+h} - f_h \rangle \right) \rho(\cdot, dx) \right) d\lambda$$

does not depend on $a, h > 0$ for any choice of $n \in \mathbb{N}$, $t_1, \ldots, t_n \geq 0$ and $\theta_1, \ldots, \theta_n \in \mathbb{R}$.

Notice that $X(0) = 0$, if and only if $mf_0 = 0$ a.e. $(\lambda)$ and $xf_0(s) = 0$ a.e. $(F)$.

EXAMPLE 1. Define a Lévy measure $\mu$ on $\mathbb{R}$ with $\int_{-1}^1 |x| \, d\mu(x) < \infty$, by $\mu((-\infty, -x)) = r(-x)$ and $\mu((x, \infty)) = r(x)$ for $x \geq 0$ $[\mu(\{0\}) = 0]$, for a nonnegative $r \in \mathbb{L}^0(\mathbb{R}) \cap \mathbb{L}^1([-1, 1])$ that is monotone and vanishes at infinity on both half-lines.

Pick an $H > 0$. Let $M$ be an i.d. random measure on $(0, \infty)$ (equipped with the Borel $\sigma$-algebra), with Lebesgue control measure, and with

$$\rho(s, B) = H s^{-1} \mu(s^{-H} B) \quad \text{and} \quad m(s) = \int_{\mathbb{R}} \kappa(x) \rho(s, dx).$$

Pick $f \in \mathbb{L}^0((0, \infty))$ satisfying (2.4). Consider the i.d. process $X$ in (2.7), where

$$(2.10) \qquad f_t(s) = \begin{cases} f(s/t), & \text{if } t > 0, \\ 0, & \text{if } t = 0, \end{cases} \qquad \text{for } s > 0.$$

This process $X$ is $H$-s.s., since the integral (2.9) evaluates to

$$\int_0^\infty \int_{\mathbb{R}} \left( \exp\left( i x a^{-H} \sum_{j=1}^n \theta_j f(s/(a t_j)) \right) - 1 \right) \frac{H \mu(s^{-H} dx)}{s} \, ds$$

$$= \int_0^\infty \int_{\mathbb{R}} \left( \exp\left( i \tilde{x} \sum_{j=1}^n \theta_j f(\tilde{s}/t_j) \right) - 1 \right) \frac{H \mu(\tilde{s}^{-H} d\tilde{x})}{\tilde{s}} \, d\tilde{s}.$$

Moreover, we get that $X$ is $\mathbf{P}$-continuous, from the fact that

$$\mathbf{E}\{e^{i\theta(X(t+h) - X(t))}\}$$

$$= \exp\left\{ \int_0^\infty \int_{\mathbb{R}} (e^{i\theta x(f(s/(t+h)) - f(s/t))} - 1) \frac{H\mu(s^{-H} dx)}{s} \, ds \right\}.$$

If $f = \mathbf{1}_{(0,1]}$, then $X$ has independent increments, so that $\mathbf{P}$-continuity and separability give local boundedness [e.g., Sato (1999), Theorem 11.5].

We conclude this section by stating a special case of Rosiński and Samorodnitsky [(1993), Theorem 2.1], that is sufficient for our needs, for easy reference.

THEOREM A.  *Let $\{X(t)\}_{t \in T}$ be an i.d. stochastic process with no Gaussian component, and with Lévy measure $\nu$ given by (2.1). Assume that the parameter space $T$ is countable, and that*

$$\mathbf{P}\left\{ \sup_{t \in T} |Z(t)| < \infty \right\} = 1.$$

*If the function*

$$H(u) \equiv 1 \wedge \nu\left( \left\{ y \in \mathbb{R}^T : \sup_{t \in T} y(t) > u \right\} \right)$$

*is subexponential, then we have*

$$\lim_{u \to \infty} \frac{1}{H(u)} \mathbf{P}\left\{ \sup_{t \in T} Z(t) > u \right\} = 1.$$

2.4. *Representation of $H$-s.s. $\alpha$-stable processes.* Let $X$ be strictly $\alpha$-stable $H$-s.s., with $\alpha \in (0, 2)$. In the case $\alpha = 1$, assume in addition that $X$ is $S\alpha S$.

Let $w_0 \in \mathbb{L}^0(S)$ be positive, and pick a constant $\beta \in [-1, 1]$ ($\beta = 0$ if $\alpha = 1$). Let $M$ be an i.d. random measure (see Section 2.3), with control measure $\lambda$, and with

$$\rho(s, dx) = w_0(s)((1 - \beta)\mathbf{1}_{(-\infty, 0)}(x) + (1 + \beta)\mathbf{1}_{(0, \infty)}(x))\frac{dx}{(2\alpha)|x|^{\alpha+1}},$$

$$m(s) = \begin{cases} \displaystyle\int_{\mathbb{R}} \kappa(x)\rho(s, dx), & \text{if } \alpha < 1, \\ 0, & \text{if } \alpha = 1, \\ \displaystyle\int_{\mathbb{R}} [\kappa(x) - x]\rho(s, dx), & \text{if } \alpha > 1. \end{cases}$$

We say that $M$ is a strictly $\alpha$-*stable random measure.* [It is an exercise to deduce from (2.2) that $M(A)$ is strictly $\alpha$-stable ($S\alpha S$ if $\alpha = 1$), for $A \in \mathfrak{S}_0$.]

There exists $\{f_t\}_{t \geq 0} \subseteq \mathbb{L}^\alpha(S, w_0\lambda)$ [which is what (2.4) reduces to here], such that $X$ satisfies (2.7), for some $\beta \in [-1, 1]$ (e.g., $\beta = \pm 1$ works if $\alpha \neq 1$).

Now a process given by (2.7) is strictly $\alpha$-stable. Denoting $x^{\langle\alpha\rangle} = |x|^\alpha \operatorname{sign}(x)$, (2.5) shows that $X$ is $H$-s.s., if and only if $f_0 = 0$ a.e. ($\lambda$), and the following integrals do not depend on $a > 0$ for any choice of $n \in \mathbb{N}$, $t_1, \ldots, t_n \geq 0$ and $\theta_1, \ldots, \theta_n \in \mathbb{R}$:

$$\int_S \left| a^{-H} \sum_{j=1}^n \theta_j f_{at_j} \right|^\alpha w_0 \, d\lambda \quad \text{and} \quad \beta \int_S \left( a^{-H} \sum_{j=1}^n \theta_j f_{at_j} \right)^{\langle\alpha\rangle} w_0 \, d\lambda.$$

Further, $X$ is $H$-s.s.s.i., if and only if $f_0 = 0$ a.e. ($\lambda$), and the following integrals do not depend on $a, h > 0$ for any choice of $n \in \mathbb{N}$, $t_1, \ldots, t_n \geq 0$ and $\theta_1, \ldots, \theta_n \in \mathbb{R}$:

$$\int_S \left| a^{-H} \sum_{j=1}^n \theta_j(f_{at_j+h} - f_h) \right|^\alpha w_0 \, d\lambda$$

and

$$\beta \int_S \left( a^{-H} \sum_{j=1}^n \theta_j(f_{at_j+h} - f_h) \right)^{\langle\alpha\rangle} w_0 \, d\lambda.$$

REMARK. Much is known about the class of $H$-s.s.s.i. $\alpha$-stable processes, that is very rich for $\alpha < 2$, unlike the Gaussian case. See, for example, Samorodnitsky and Taqqu (1994), Surgailis, Rosiński, Mandrekar and Cambanis (1998), Burnecki, Rosiński and Weron (1998) and Pipiras and Taqqu (2002a, b).

For $H \in (1/\alpha, 1]$ with $\alpha > 1$, and for $H = 1/\alpha > 1$, it is known that (separable) $H$-s.s.s.i. $\alpha$-stable processes are locally bounded. For other values of $H$ and $\alpha$, local boundedness is not determined by $H$ and $\alpha$, and there exist both locally bounded and unbounded processes. Precise conditions for local boundedness are known for $\alpha < 1$. See Kôno and Maejima (1991) and Samorodnitsky and Taqqu (1990, 1994).

**3. Overload and the Piterbarg properties.** Here we first study the probability for overload $\mathbf{P}\{Y(t) > u\}$, and then the Piterbarg properties (1.2).

The next assumptions limit the effect of the left tail of $X$ on the right tail of $Y$:

$$(3.1) \quad \limsup_{u \to \infty} \frac{\mathbf{P}\{X(1) < -\varepsilon (t(u) u^{-1/\gamma})^{-H} u^{1-H/\gamma}\}}{R(u^{1-H/\gamma})} < \infty \qquad \text{for all } \varepsilon > 0,$$

$$(3.2) \quad \lim_{u \to \infty} \frac{\mathbf{P}\{X(1) < -\varepsilon (t(u) u^{-1/\gamma})^{-H} v(u^{1-H/\gamma})\}}{R(u^{1-H/\gamma})} = 0 \qquad \text{for some } \varepsilon > 0.$$

Here (3.1) is used together with the growth condition indicated in the Introduction

$$(3.3) \quad \limsup_{u \to \infty} \frac{t(u)}{u^{1/\gamma}} < \infty.$$

In (3.2), $v : \mathbb{R} \to (0, \infty)$ is a suitably selected function, with the following properties:

$$(3.4) \quad \limsup_{u \to \infty} \frac{v(u)}{u} < \infty \quad \text{and} \quad \lim_{u \to \infty} \frac{t(u)}{u^{1/\gamma} (v(u^{1-H/\gamma})/u^{1-H/\gamma})^{1/(H \wedge 1)}} = 0.$$

In the hypothesis of Theorem 1, (3.1) follows from (3.3) [by (3.11) and $R \in OR$], while (3.2) follows from (3.4) [by (3.20) and $R \in PD$], provided that

$$(3.5) \quad \limsup_{u \to \infty} \frac{\mathbf{P}\{X(1) < -u\}}{\mathbf{P}\{\sup_{s \in [0,1]} X(s) > u\}} < \infty.$$

Notice that (3.5) holds if, for example, $X(1)$ is symmetric or nonnegative.

THEOREM 1. *Let $X$ be $H$-s.s. and i.d. with no Gaussian component, satisfying Condition $\mathfrak{X}$. Consider the process $Y$, given by (1.1), with $c > 0$ and $\gamma > H$ constants. Suppose that the function $R$, given by (1.3) [with $\nu$ given by (2.1)], belongs to $\mathcal{S} \cap PD$. Then $Y(t) < \infty$ a.s. for each $t \geq 0$, and*

$$(3.6) \quad \lim_{u \to \infty} \frac{\mathbf{P}\{Y(0) > u\}}{R(u^{1-H/\gamma})} = 1.$$

(i) *If $X$ is s.i., then for $t(u) \geq 0$,*

$$\lim_{u \to \infty} \frac{\mathbf{P}\{Y(t(u)) > u\}}{R(u^{1-H/\gamma})} = 1. \tag{3.7}$$

(ii) *If $R \in OR$ and $t(u) \geq 0$ satisfies (3.1) and (3.3), then we have*

$$0 < \liminf_{u \to \infty} \frac{\mathbf{P}\{Y(t(u)) > u\}}{R(u^{1-H/\gamma})} \leq \limsup_{u \to \infty} \frac{\mathbf{P}\{Y(t(u)) > u\}}{R(u^{1-H/\gamma})} < \infty. \tag{3.8}$$

(iii) *If there exists a function $v$ satisfying (3.4), such that*

$$\lim_{u \to \infty} \frac{R(u - v(u))}{R(u)} = 1, \tag{3.9}$$

*and such that (3.2) holds, then (3.7) holds.*

PROOF. Let $R \in \mathcal{S} \cap PD$, and denote $\hat{c} = 2^{\gamma}/c$. For every $t \geq 0$, we have, for $u > 0$,

$$\mathbf{P}\{Y(t) > u\} \leq \mathbf{P}\left\{X(t) < -\frac{u}{2}\right\}$$
$$+ \mathbf{P}\left\{\sup_{t \leq s < 2t} X(s) > \frac{u}{2}\right\} + \mathbf{P}\left\{\sup_{s \geq 2t}(X(s) - s^{\gamma}/\hat{c}) > \frac{u}{2}\right\}.$$

Since $X$ is locally bounded, the first two terms on the right go to zero as $u \to \infty$. Furthermore, we can bound from above the third term by

$$\mathbf{P}\left\{\sup_{0 \leq s \leq u^{1/\gamma}} X(s) > \frac{u}{2}\right\} + \sum_{j=0}^{\infty} \mathbf{P}\left\{\sup_{2^j u^{1/\gamma} \leq s \leq 2^{j+1} u^{1/\gamma}} X(s) > \frac{2^{\gamma j} u}{\hat{c}}\right\}$$
$$= \mathbf{P}\left\{\sup_{s \in [0,1]} X(s) > \frac{u^{1-H/\gamma}}{2}\right\} \tag{3.10}$$
$$+ \sum_{j=0}^{\infty} \mathbf{P}\left\{\sup_{s \in [1/2,1]} X(s) > \frac{2^{(\gamma-H)j} u^{1-H/\gamma}}{\hat{c}}\right\}$$
$$\leq 2 \sum_{j=0}^{\infty} \mathbf{P}\left\{\sup_{s \in [0,1]} \frac{X(s)}{1 + cs^{\gamma}} > \frac{2^{(\gamma-H)j} u^{1-H/\gamma}}{(1+c)(2 \vee \hat{c})}\right\}.$$

Here $R \in \mathcal{S}$, together with Theorem A, gives

$$\limsup_{u \to \infty} \frac{1}{R(u)} \mathbf{P}\left\{\sup_{s \in [0,1]} X(s) > (1+c)u\right\} \tag{3.11}$$
$$\leq \limsup_{u \to \infty} \frac{1}{R(u)} \mathbf{P}\left\{\sup_{s \in [0,1]} \frac{X(s)}{1 + cs^{\gamma}} > u\right\} \leq 1.$$

Here and in future applications of Theorem A, we use that the process under consideration is separable and $\mathbf{P}$-continuous. Hence it is enough to consider suprema over any countable dense subset of the parameter space of that process [e.g., Samorodnitsky and Taqqu (1994), Exercise 9.3], which we take to be the rational numbers in the interior of the parameter space (when that parameter space is an interval).

Our $R \in PD$ has a so-called upper Matuszewska index $\mathfrak{a} < 0$ [e.g., Bingham, Goldie and Teugel (1987), page 71] such that, given $-a \in (\mathfrak{a}, 0)$ and $\lambda_0 > 0$,

$$\frac{R(\lambda u)}{R(u)} \leq C\lambda^{-a} \qquad \text{for } \lambda \geq \lambda_0 \text{ and } u \text{ large enough,}$$

for some $C > 0$. Hence the right-hand side of (3.10) is at most

$$2\sum_{j=0}^{\infty} 2R\left(\frac{2^{(\gamma - H)j} u^{1 - H/\gamma}}{(1 + c)(2 \vee \hat{c})}\right)$$

$$\leq 4C(1 + c)^a (2 \vee \hat{c})^a \left(\sum_{j=0}^{\infty} 2^{-a(\gamma - H)j}\right) R(u^{1 - H/\gamma})$$

for $u$ large enough. This proves the fact that $Y(t) < \infty$ a.s.   $\square$

Further, by self-similarity and Theorem A [cf. (3.11)],

$$(3.12) \qquad \begin{aligned} \mathbf{P}\{Y(0) > u\} &= \mathbf{P}\left\{\sup_{s \geq 0} \frac{X(s)}{1 + cs^{\gamma}} > u^{1 - H/\gamma}\right\} \\ &\sim R(u^{1 - H/\gamma}) \qquad \text{as } u \to \infty. \end{aligned}$$

PROOF OF (i).   For $X$ s.i., $Y$ is stationary, and so (3.7) is the same thing as (3.6).   $\square$

PROOF OF (ii).   By (3.3), we have, for some $\theta \in (0, 1)$, for all $u$ large enough,

$$(3.13) \quad \theta \leq \inf_{s \geq tu^{-1/\gamma}} \frac{1 + c(s - tu^{-1/\gamma})^{\gamma}}{1 + cs^{\gamma}} \leq \sup_{s \geq tu^{-1/\gamma}} \frac{1 + c(s - tu^{-1/\gamma})^{\gamma}}{1 + cs^{\gamma}} \leq 1$$

[where $t = t(u)$]. Using self-similarity, we therefore obtain, for $u$ large enough,

$$(3.14) \qquad \begin{aligned} \mathbf{P}\{Y(t) > u\} &= \mathbf{P}\left\{\sup_{s \geq tu^{-1/\gamma}} \frac{X(s) - X(tu^{-1/\gamma})}{1 + c(s - tu^{-1/\gamma})^{\gamma}} > u^{1 - H/\gamma}\right\} \\ &\begin{cases} \geq \mathbf{P}\left\{\sup_{s \geq tu^{-1/\gamma}} \frac{X(s) - X(tu^{-1/\gamma})}{1 + cs^{\gamma}} > u^{1 - H/\gamma}\right\}, \\ \leq \mathbf{P}\left\{\sup_{s \geq 0} \frac{X(s) - X(tu^{-1/\gamma})}{1 + cs^{\gamma}} > \theta u^{1 - H/\gamma}\right\}. \end{cases} \end{aligned}$$

Notice that, denoting

$$\eta_1(u) = (tu^{-1/\gamma})^{-H}(1 + c(Ktu^{-1/\gamma})^\gamma),$$

$$\eta_2(u) = (tu^{-1/\gamma})^{-H}(1 \vee (Ktu^{-1/\gamma}))^\gamma,$$

for a constant $K \geq 1$, we obtain

$$\mathbf{P}\left\{\sup_{s \geq tu^{-1/\gamma}} \frac{X(s) - X(tu^{-1/\gamma})}{1 + cs^\gamma} > u^{1-H/\gamma}\right\}$$

$$\geq \mathbf{P}\left\{\sup_{s \geq Ktu^{-1/\gamma}} \frac{X(s) - X(tu^{-1/\gamma})}{1 + cs^\gamma} > u^{1-H/\gamma}\right\}$$

$$\geq \mathbf{P}\left\{\sup_{s \geq Ktu^{-1/\gamma}} \frac{X(s)}{1 + cs^\gamma} > 2u^{1-H/\gamma}\right\}$$

$$\quad - \mathbf{P}\{X(tu^{-1/\gamma}) > (1 + c(Ktu^{-1/\gamma})^\gamma)u^{1-H/\gamma}\}$$

$$= \mathbf{P}\left\{\sup_{s \geq 1} \frac{X(s)}{1 + c(Ktu^{-1/\gamma})^\gamma s^\gamma} > 2(Ktu^{-1/\gamma})^{-H}u^{1-H/\gamma}\right\}$$

(3.15)

$$\quad - \mathbf{P}\{X(1) > \eta_1(u)u^{1-H/\gamma}\}$$

$$\geq \mathbf{P}\left\{\sup_{s \geq 1} \frac{X(s)}{1 + cs^\gamma} > \frac{2\eta_2(u)u^{1-H/\gamma}}{K^H}\right\}$$

$$\quad - \mathbf{P}\left\{\sup_{s \geq 0} \frac{X(s)}{1 + cs^\gamma} > \frac{\eta_1(u)u^{1-H/\gamma}}{1 + c}\right\}.$$

Here we have, picking a constant $L \geq 1$,

$$\mathbf{P}\left\{\sup_{s \geq 0} \frac{X(s)}{1 + cs^\gamma} > u\right\}$$

$$\leq \mathbf{P}\left\{\sup_{s \geq L^{-1}} \frac{X(s)}{1 + cs^\gamma} > u\right\} + \mathbf{P}\left\{\sup_{0 \leq s \leq L^{-1}} \frac{X(s)}{1 + cs^\gamma} > u\right\}$$

$$= \mathbf{P}\left\{\sup_{s \geq 1} \frac{X(s/L)}{1 + c(s/L)^\gamma} > u\right\} + \mathbf{P}\left\{\sup_{0 \leq s \leq 1} \frac{X(s/L)}{1 + cs^\gamma} \frac{1 + cs^\gamma}{1 + c(s/L)^\gamma} > u\right\}$$

$$\leq \mathbf{P}\left\{\sup_{s \geq 1} \frac{X(s)}{1 + cs^\gamma} > \frac{u}{L^{\gamma - H}}\right\} + \mathbf{P}\left\{\sup_{s \geq 0} \frac{X(s)}{1 + cs^\gamma} > \frac{L^H u}{1 + c}\right\}.$$

It follows from (3.12) and the fact that $R \in PD$ that, if $L$ is large enough, then

$$\limsup_{u \to \infty} \mathbf{P}\left\{\sup_{s \geq 0} \frac{X(s)}{1 + cs^\gamma} > \frac{L^H u}{1 + c}\right\} \Big/ \mathbf{P}\left\{\sup_{s \geq 0} \frac{X(s)}{1 + cs^\gamma} > u\right\} < \frac{1}{2}.$$

Fixing $L$ such that this relation holds, we get immediately

$$\mathbf{P}\left\{\sup_{s\geq 1}\frac{X(s)}{1+cs^\gamma} > \frac{u}{L^{\gamma-H}}\right\} \geq \frac{1}{2}\mathbf{P}\left\{\sup_{s\geq 0}\frac{X(s)}{1+cs^\gamma} > u\right\}$$

for $u$ large enough. Therefore, by (3.12),

$$(3.16) \qquad \liminf_{u\to\infty}\frac{1}{R(L^{\gamma-H}u)}\mathbf{P}\left\{\sup_{s\geq 1}\frac{X(s)}{1+cs^\gamma} > u\right\} \geq \frac{1}{2}.$$

Since $R \in PD$, and $\limsup_{u\to\infty}\eta_2(u)/\eta_1(u) < \infty$, we get from (3.12) and (3.16),

$$(3.17) \qquad \begin{aligned} &\liminf_{u\to\infty}\mathbf{P}\left\{\sup_{s\geq 1}\frac{X(s)}{1+cs^\gamma} > \frac{2\eta_2(u)u^{1-H/\gamma}}{K^H}\right\} \\ &\qquad \times \mathbf{P}\left\{\sup_{s\geq 0}\frac{X(s)}{1+cs^\gamma} > \frac{\eta_1(u)u^{1-H/\gamma}}{1+c}\right\}^{-1} \\ &\geq \frac{1}{2}\liminf_{u\to\infty}R\left(\frac{2L^{\gamma-H}\eta_2(u)u^{1-H/\gamma}}{K^H}\right)\Big/R\left(\frac{\eta_1(u)u^{1-H/\gamma}}{1+c}\right) \\ &\geq 2 \end{aligned}$$

for all $K$ large enough [where $\eta_i(u)u^{1-H/\gamma} \to \infty$ for $i = 1, 2$]. Fixing $K \geq 1$ such that (3.17) holds, we may apply (3.16) and (3.17) on the last row of (3.15), to get

$$\begin{aligned} &\liminf_{u\to\infty}\frac{1}{R(u^{1-H/\gamma})}\mathbf{P}\left\{\sup_{s\geq tu^{-1/\gamma}}\frac{X(s)-X(tu^{-1/\gamma})}{1+cs^\gamma} > u^{1-H/\gamma}\right\} \\ &\geq \frac{1}{2}\liminf_{u\to\infty}\frac{1}{R(u^{1-H/\gamma})}\mathbf{P}\left\{\sup_{s\geq 1}\frac{X(s)}{1+cs^\gamma} > \frac{2\eta_2(u)u^{1-H/\gamma}}{K^H}\right\} \\ &\geq \frac{1}{4}\liminf_{u\to\infty}\frac{R(2K^{-H}L^{\gamma-H}\eta_2(u)u^{1-H/\gamma})}{R(u^{1-H/\gamma})} \\ &> 0, \end{aligned}$$

using $R \in OR$ for the last inequality. By (3.14), this gives the lower bound in (3.8).

The corresponding upper bound in (3.8) follows from (3.14). This is so because

$$\begin{aligned} &\limsup_{u\to\infty}\frac{1}{R(u^{1-H/\gamma})}\mathbf{P}\left\{\sup_{s\geq 0}\frac{X(s)-X(tu^{-1/\gamma})}{1+cs^\gamma} > \theta u^{1-H/\gamma}\right\} \\ &\leq \limsup_{u\to\infty}\frac{1}{R(u^{1-H/\gamma})}\mathbf{P}\left\{\sup_{s\geq 0}\frac{X(s)}{1+cs^\gamma} > \frac{\theta}{2}u^{1-H/\gamma}\right\} \\ &\qquad + \limsup_{u\to\infty}\frac{\mathbf{P}\{-X(tu^{-1/\gamma}) > (\theta/2)u^{1-H/\gamma}\}}{R(u^{1-H/\gamma})}, \end{aligned}$$

which is finite, by (3.1) and (3.12), together with self-similarity and $R \in OR$. $\square$

PROOF OF (iii). We have

$$
\begin{aligned}
(3.18) \qquad 1 - O((tu^{-1/\gamma})^{\gamma \wedge 1}) &\leq \inf_{s \geq tu^{-1/\gamma}} \frac{1 + c(s - tu^{-1/\gamma})^{\gamma}}{1 + cs^{\gamma}} \\
&\leq \sup_{s \geq tu^{-1/\gamma}} \frac{1 + c(s - tu^{-1/\gamma})^{\gamma}}{1 + cs^{\gamma}} \\
&\leq 1
\end{aligned}
$$

as $u \to \infty$. Here (3.4) shows that, with obvious notation,

$$
\begin{aligned}
O((tu^{-1/\gamma})^{\gamma \wedge 1}) &\leq o\left((v(u^{1-H/\gamma})/u^{1-H/\gamma})^{(\gamma \wedge 1)/(H \wedge 1)}\right) \\
&\leq o(v(u^{1-H/\gamma})/u^{1-H/\gamma}).
\end{aligned}
$$

This gives us the following version of (3.14), that for $u$ large enough:

$$
\mathbf{P}\{Y(t) > u\}
\begin{cases}
\geq \mathbf{P}\left\{\sup_{s \geq tu^{-1/\gamma}} \dfrac{X(s) - X(tu^{-1/\gamma})}{1 + cs^{\gamma}} > u^{1-H/\gamma}\right\}, \\
\leq \mathbf{P}\left\{\sup_{s \geq 0} \dfrac{X(s) - X(tu^{-1/\gamma})}{1 + cs^{\gamma}} > u^{1-H/\gamma} - \varepsilon v(u^{1-H/\gamma})\right\}.
\end{cases}
$$

(3.19)

To bound the ratio in (3.7) from below, use self-similarity, (3.4) and (3.12), to get

$$
\begin{aligned}
(3.20) \qquad &\mathbf{P}\left\{\sup_{s \in [0, tu^{-1/\gamma}]} X(s) > \varepsilon v(u^{1-H/\gamma})\right\} \\
&\leq \mathbf{P}\left\{\sup_{s \geq 0} \frac{X(s)}{1 + cs^{\gamma}} > \frac{\varepsilon v(u^{1-H/\gamma})}{(1 + c)(tu^{-1/\gamma})^{H}}\right\} \\
&\leq \mathbf{P}\left\{\sup_{s \geq 0} \frac{X(s)}{1 + cs^{\gamma}} > Ku^{1-H/\gamma}\right\} \\
&\sim R(Ku^{1-H/\gamma}) \qquad \text{as } u \to \infty,
\end{aligned}
$$

for any constant $K \geq 1$. Hence (3.19), together with (3.9) and (3.12), give that

$$
\begin{aligned}
&\liminf_{u \to \infty} \frac{\mathbf{P}\{Y(t) > u\}}{R(u^{1-H/\gamma})} \\
&\qquad \geq \liminf_{u \to \infty} \frac{1}{R(u^{1-H/\gamma})} \mathbf{P}\left\{\sup_{s \geq 0} \frac{X(s)}{1 + cs^{\gamma}} > u^{1-H/\gamma} - 2\varepsilon v(u^{1-H/\gamma})\right\} \\
&\qquad\quad - \limsup_{u \to \infty} \frac{\mathbf{P}\{X(tu^{-1/\gamma}) > \varepsilon v(u^{1-H/\gamma})\}}{R(u^{1-H/\gamma})}
\end{aligned}
$$

$$- \limsup_{u \to \infty} \frac{1}{R(u^{1-H/\gamma})} \mathbf{P} \left\{ \sup_{s \in [0, tu^{-1/\gamma})} X(s) > \varepsilon v(u^{1-H/\gamma}) \right\}$$

$$\geq 1 - 2 \limsup_{u \to \infty} \frac{R(Ku^{1-H/\gamma})}{R(u^{1-H/\gamma})} \to 1 \qquad \text{as } K \to \infty,$$

since $R \in PD$. Of course, this establishes that

$$\liminf_{u \to \infty} \frac{\mathbf{P}\{Y(t) > u\}}{R(u^{1-H/\gamma})} \geq 1.$$

On the other hand, since (3.9) and monotonicity of $R$ give $R(u - \lambda v(u)) \sim R(u)$ for any $\lambda \in \mathbb{R}$, (3.19) together with (3.2), (3.9) and (3.12), show that

$$\limsup_{u \to \infty} \frac{\mathbf{P}\{Y(t) > u\}}{R(u^{1-H/\gamma})}$$

$$\leq \limsup_{u \to \infty} \frac{1}{R(u^{1-H/\gamma})} \mathbf{P} \left\{ \sup_{s \geq 0} \frac{X(s)}{1 + cs^\gamma} > u^{1-H/\gamma} - 2\varepsilon v(u^{1-H/\gamma}) \right\}$$

$$\quad + \limsup_{u \to \infty} \frac{\mathbf{P}\{-X(1) > \varepsilon(tu^{-1/\gamma})^{-H} v(u^{1-H/\gamma})\}}{R(u^{1-H/\gamma})}$$

$$= 1 + 0. \hspace{6cm} \square$$

To establish the Piterbarg property, we make use of the following assumptions:

$$(3.21) \qquad \limsup_{u \to \infty} \frac{\mathbf{P}\{\inf_{s \in [0,1]} X(s) < -\varepsilon(t(u)u^{-1/\gamma})^{-H} u^{1-H/\gamma}\}}{R(u^{1-H/\gamma})} < \infty$$

$$\text{for all } \varepsilon > 0,$$

$$(3.22) \qquad \lim_{u \to \infty} \frac{\mathbf{P}\{\inf_{s \in [0,1]} X(s) < -\varepsilon(t(u)u^{-1/\gamma})^{-H} v(u^{1-H/\gamma})\}}{R(u^{1-H/\gamma})} = 0$$

$$\text{for some } \varepsilon > 0.$$

Assumption (3.21) will be used together with the growth condition (3.3), while in assumption (3.22), $v$ is a suitably selected function that satisfies (3.4).

In the hypothesis of Theorem 2, (3.21) follows from (3.3) [by (3.11) and $R \in OR$], while (3.22) follows from (3.4) [by (3.20) and monotonicity of $R$], when

$$(3.23) \qquad \limsup_{u \to \infty} \frac{\mathbf{P}\{\inf_{s \in [0,1]} X(s) < -u\}}{\mathbf{P}\{\sup_{s \in [0,1]} X(s) > u\}} < \infty.$$

Clearly, (3.23) holds for $X$ symmetric or nonnegative. Otherwise, (3.21) and (3.22) could possibly be verified by Theorem A, for $-\inf_{s \in [0,1]} X(s)$ subexponential, or by Albin [(1998), Theorem 3 and Sections 8 and 9].

THEOREM 2. *Let $X$ be $H$-s.s. and i.d. with no Gaussian component, satisfying Condition $\mathfrak{X}$. Consider the process $Y$, given by (1.1), together with the function $R$, given by (1.3), where $c > 0$ and $\gamma > H$ are constants.*

(i) *Let $R \in \mathcal{S} \cap OR \cap PD$. If (3.3) and (3.21) hold, we have, for $0 \leq \hat{t}(u) \leq t(u)$,*

$$1 \leq \liminf_{u \to \infty} \frac{\mathbf{P}\{\sup_{s \in [0, t(u)]} Y(s) > u\}}{\mathbf{P}\{Y(\hat{t}(u)) > u\}}$$

$$\leq \limsup_{u \to \infty} \frac{\mathbf{P}\{\sup_{s \in [0, t(u)]} Y(s) > u\}}{\mathbf{P}\{Y(\hat{t}(u)) > u\}} < \infty.$$

(ii) *Let $R \in \mathcal{S} \cap PD$. Take $t(u)$ and $v$ such that (3.4), (3.9) and (3.22) hold. The process $Y$ has the strong Piterbarg property (1.2c).*

PROOF OF (i). It is enough to show the upper bound. Using (3.13) and (3.14), we get

$$\mathbf{P}\left\{\sup_{s \in [0, t]} Y(s) > u\right\}$$
$$= \mathbf{P}\left\{\sup_{0 \leq r \leq tu^{-1/\gamma}} \sup_{s \geq r} \frac{X(s) - X(r)}{1 + c(s - r)^\gamma} > u^{1 - H/\gamma}\right\}$$
$$\leq \mathbf{P}\left\{\sup_{s \geq 0} \frac{X(s)}{1 + cs^\gamma} > \frac{\theta u^{1 - H/\gamma}}{2}\right\} + \mathbf{P}\left\{\sup_{0 \leq r \leq tu^{-1/\gamma}} -X(r) > \frac{\theta u^{1 - H/\gamma}}{2}\right\}$$
(3.24)

for some $\theta \in (0, 1)$. Therefore, self-similarity, (3.12), (3.21) and $R \in OR$, give

$$\limsup_{u \to \infty} \frac{\mathbf{P}\{\sup_{s \in [0, t]} Y(s) > u\}}{R(u^{1 - H/\gamma})} < \infty.$$

Now the upper bound desired follows from (3.8) [notice that (3.21) implies (3.1)].

$\square$

PROOF OF (ii). Using (3.18) together with (3.4), as in the last paragraph of the proof of part (iii) of Theorem 1, we may readily modify the estimate (3.24) to obtain

$$\mathbf{P}\left\{\sup_{s \in [0, t]} Y(s) > u\right\}$$
$$\leq \mathbf{P}\left\{\sup_{s \geq 0} \frac{X(s)}{1 + cs^\gamma} > u^{1 - H/\gamma} - \varepsilon v(u^{1 - H/\gamma})\right\}$$
$$+ \mathbf{P}\left\{\sup_{0 \leq r \leq tu^{-1/\gamma}} -X(r) > \varepsilon v(u^{1 - H/\gamma})\right\}.$$

By application of (3.9) together with (3.12) and (3.22), this shows that

$$\limsup_{u \to \infty} \frac{\mathbf{P}\{\sup_{s \in [0,t]} Y(s) > u\}}{R(u^{1-H/\gamma})} \leq 1.$$

On the other hand,

$$\mathbf{P}\left\{\inf_{s \in [0,t]} Y(s) > u\right\}$$

$$= \mathbf{P}\left\{\inf_{0 \leq r \leq tu^{-1/\gamma}} \sup_{s \geq r} \frac{X(s) - X(r)}{1 + c(s-r)^\gamma} > u^{1-H/\gamma}\right\}$$

$$\geq \mathbf{P}\left\{\inf_{0 \leq r \leq tu^{-1/\gamma}} \left(\sup_{s \geq r} \frac{X(s)}{1 + c(s-r)^\gamma} - X(r)_+\right) > u^{1-H/\gamma}\right\}$$

$$\geq \mathbf{P}\left\{\inf_{0 \leq r \leq tu^{-1/\gamma}} \sup_{s \geq r} \frac{X(s)}{1 + cs^\gamma} - \sup_{0 \leq r \leq tu^{-1/\gamma}} X(r)_+ > u^{1-H/\gamma}\right\}$$

$$= \mathbf{P}\left\{\sup_{s \geq tu^{-1/\gamma}} \frac{X(s)}{1 + cs^\gamma} - \sup_{0 \leq r \leq tu^{-1/\gamma}} X(r)_+ > u^{1-H/\gamma}\right\},$$

and, as was established in the proof of Theorem 1, this gives us

$$\liminf_{u \to \infty} \frac{\mathbf{P}\{\inf_{s \in [0,t]} Y(s) > u\}}{R(u^{1-H/\gamma})} \geq 1.$$

Hence the strong Piterbarg property.  $\square$

Here are two easy corollaries to Theorems 1 and 2, that make use of (3.23):

COROLLARY 1. *Let $X$ be $H$-s.s. and i.d. with no Gaussian component, satisfying Condition $\mathfrak{X}$. Let the function $R$ belong to $\mathcal{S} \cap PD$, and assume that (3.23) holds.*

(i) *The Piterbarg property holds for $\gamma \leq H + (H \wedge 1)$.*

(ii) *If (3.9) holds for $v(u) = (1 \vee u)^\beta$, for some $\beta \in (0,1)$, then the Piterbarg property holds for $\gamma < H + (H \wedge 1)/(1-\beta)$. Hence, if (3.9) holds for $v(u) = (1 \vee u)^\beta$, for each $\beta \in (0,1)$, then the Piterbarg property holds for any $\gamma > H$.*

PROOF OF (i). By part (ii) of Theorem 2, together with an inspection of (3.4), it is enough to exhibit a positive function $v$, with $\lim_{u \to \infty} v(u) = \infty$ and $\limsup_{u \to \infty} v(u)/u < \infty$, that satisfies (3.9). This is easy: Let $b_0 = 0$, $b_1 = 1$ and

$$b_{i+1} = \inf\left\{u \geq \max(b_i, 2i) : \inf_{x \geq u} \frac{R(x)}{R(x-i)} \geq 1 - \frac{1}{i}\right\} \qquad \text{for } i \geq 1.$$

Since $R \in \mathcal{S} \subseteq \mathcal{L}$ (see Section 2.1), this is an increasing to infinity sequence of finite nonnegative numbers, and we may now choose

$$v(u) = i \quad \text{if } u \in [b_i, b_{i+1}) \text{ for } i \geq 1 \quad \text{and} \quad v(u) = 1 \text{ for } u < 1. \qquad \square$$

PROOF OF (ii). Once again, the result follows from part (ii) of Theorem 2, by means of checking that (3.4) holds for $\gamma < H + (H \wedge 1)/(1 - \beta)$, when $v(u) = (u \vee 1)^\beta$. $\square$

In Part II of Example 5, below we see that the Piterbarg property may be absent, when, in the notation of part (ii) of Corollary 1, $\gamma \geq H + 1/(1 - \beta)$.

COROLLARY 2. *Let $X$ be $H$-s.s. and i.d. with no Gaussian component, satisfying Condition $\mathfrak{X}$. Let the function $R$ belong to $\mathcal{S} \cap PD$, and assume that (3.23) holds.*

(i) *The strong Piterbarg property holds for $t(u) \geq 0$ such that*

$$\limsup_{u \to \infty} \frac{t(u)}{u^{1/\gamma - (1 - H/\gamma)/(H \wedge 1)}} < \infty.$$

(ii) *If (3.9) holds for $v(u) = (1 \vee u)^\beta$, for some $\beta \in (0, 1)$, then the strong Piterbarg property holds for $t(u) \geq 0$ such that*

$$(3.25) \qquad \lim_{u \to \infty} \frac{t(u)}{u^{1/\gamma - (1 - \beta)(1 - H/\gamma)/(H \wedge 1)}} = 0.$$

(iii) *If (3.9) holds for every $v(u) = o(u)$, then the strong Piterbarg property holds for each $t(u) = o(u^{1/\gamma})$.*

Corollary 2 is proved in the same way as Corollary 1.

Notice that the intervals with length $t(u)$ as in (3.25), for which the strong Piterbarg property holds, do in fact shrink with $u$, unless $\gamma < H + (H \wedge 1)/(1 - \beta)$.

**4. Discussion and examples.** Here we present points of view on the results of Section 3. Examples are given, where the input process $X$ is represented in the form of a stochastic integral, with respect to an i.d. random measure.

We first exemplify that the storage process does not have to be finite valued, in general. (Obviously, this does not happen under the assumptions of our results.)

EXAMPLE 2. For standard Brownian $\{B(t)\}_{t \geq 0}$ motion, and a nondecreasing function $f : (0, \infty) \to (0, \infty)$, by the Kolmogorov–Petrowski integral test,

$$\{B(t) \leq \sqrt{2t} f(t) \text{ for } t \geq T, \text{ for some } T = T(\omega) < \infty\}$$

is a zero-one event, or in other words,

$$\mathbf{P}\{B(t) \le \sqrt{2t}\, f(t) \text{ ultimately as } t \to \infty\} = 0 \text{ or } 1,$$

with the probability being 1 if and only if

$$\int_{t_0}^{\infty} \frac{f(t)}{t} \exp\{-f(t)^2\}\, dt < \infty \qquad \text{for some } t_0 \ge 0$$

[e.g., Bingham (1986), page 436]. From this we get that

$$(4.1) \qquad \mathbf{P}\left\{\frac{B(t)^2}{2t} \le \ln\ln t + \frac{3\ln\ln\ln t}{2} \text{ ultimately as } t \to \infty\right\} = 0.$$

Consider the following $H$-s.s. process $X$:

$$(4.2) \qquad X(t) = t^H \exp\left\{\exp\left[\frac{B(t)^2}{2t}\right]\right\} \qquad \text{for } t > 0.$$

From (4.1) it follows readily that

$$\limsup_{t \to \infty} \frac{X(t)}{t^\gamma} = \infty \qquad \text{w.p.1 for } \gamma \in \mathbb{R}.$$

Hence, the storage process $Y(t)$ in (1.1), with input $X$ given by (4.2), is not finite valued for any $t \ge 0$. [Incidentally, using the fact that $\{tB(1/t)\}_{t \ge 0} \overset{d}{=} \{B(t)\}_{t \ge 0}$ together with (4.1), it can be seen that the process $X$ is not bounded at zero.]

4.1. *I.d. $H$-s.s. processes.* In this section, $X$ denotes the i.d. $H$-s.s. process given in Example 1, which is assumed to satisfy Condition $\mathfrak{X}$. Notice that, by (2.6), $X$ is nonnegative, if $f$ is nonnegative (nonpositive) and $r$ is zero on $(-\infty, 0]$ $([0, \infty))$.

Denoting, for $s > 0$,

$$\sigma_+(s)^{-1} = s^H \sup_{t \in \mathbb{Q}^+} \frac{f(s/t)^+}{1 + ct^\gamma} \quad \text{and} \quad \sigma_-(s)^{-1} = s^H \sup_{t \in \mathbb{Q}^+} \frac{f(s/t)^-}{1 + ct^\gamma},$$

we have, by (2.8), for $u$ large enough,

$$(4.3) \qquad \begin{aligned} R(u) &= H \int_0^\infty \mu(\mathbb{R} \setminus [-\sigma_-(s)u, \sigma_+(s)u]) \frac{ds}{s} \\ &= H \int_0^\infty (r(\sigma_+(s)u) + r(-\sigma_-(s)u)) \frac{ds}{s}. \end{aligned}$$

EXAMPLE 3. Let $r$ be symmetric, with $r \in RV(-\rho)$ for some $\rho > 0$. By (2.2), it follows that $X$ is symmetric. Assume that for some $\varepsilon > 0$,

$$(4.4) \qquad \int_0^\infty (\sigma_+(s)^{-(\rho+\varepsilon)} + \sigma_-(s)^{-(\rho+\varepsilon)}) \frac{ds}{s} < \infty.$$

Then by (4.3), we have

$$(4.5) \qquad R(u) \sim H \int_0^\infty (\sigma_+(s)^{-\rho} + \sigma_-(s)^{-\rho}) \frac{ds}{s} r(u) \qquad \text{as } u \to \infty.$$

Thus part (ii) of Corollary 1 shows that the Piterbarg property holds for $\gamma > H$, and part (iii) of Corollary 2 gives the strong Piterbarg property for $t(u) = o(u^{1/\gamma})$.

In Example 3, symmetry gives us (3.23), for free. Without symmetry, we may still establish the Piterbarg properties, by direct verification of (3.23) [or (3.22)].

EXAMPLE 4. Take $f$ nonnegative, and not identically zero. Assume that

$$r \in RV(-\rho) \quad \text{and} \quad r(-\cdot) \in RV(-\rho) \qquad \text{with } \limsup_{u \to \infty} \frac{r(-u)}{r(u)} < \infty,$$

for some constant $\rho > 0$. Under the condition (4.4), we see that (4.5) holds with $\sigma_-(s) = 0$. Further, as in (4.5) by Theorem A, the limit in (3.23) is

$$\left( \int_0^\infty s^{H\rho} \sup_{r \in (s,\infty) \cap \mathbb{Q}} f(r)^\rho \frac{ds}{s} \Big/ \int_0^\infty s^{H\rho} \sup_{r \in (s,\infty) \cap \mathbb{Q}} f(r)^\rho \frac{ds}{s} \right)$$
$$\times \limsup_{u \to \infty} \frac{r(-u)}{r(u)} < \infty.$$

Hence the Piterbarg properties hold in the same way as in Example 3.

EXAMPLE 5 —Part I. Pick constants $A > 0$ and $\alpha \in (0,1)$, and consider

$$(4.6) \quad r(x) = g(x)e^{-Ax^\alpha} \qquad \text{for } x \geq 0 \quad \text{and} \quad r(x) = 0 \qquad \text{for } x < 0,$$

where $g \in RV(\rho)$ with $\rho \in \mathbb{R}$. Take $f$ nonnegative, so that $X$ is nonnegative.

Let $\sigma_+(s)$ take its minimal value at a unique $\hat{s} > 0$, where $\sigma_+$ is two times continuously differentiable at $\hat{s}$, with $\sigma_+''(\hat{s}) > 0$. By Taylor expansion in (4.3), we have

$$R(u) = Hu^{-\alpha/2} \int_{-\hat{s}u^{\alpha/2}}^\infty \frac{g(\sigma_+(\hat{s} + s/u^{\alpha/2})u)}{\hat{s} + s/u^{\alpha/2}} e^{-A(\sigma_+(\hat{s}+s/u^{\alpha/2}))^\alpha u^\alpha} ds$$
$$\sim \frac{H\sigma_+(\hat{s})^\rho}{\hat{s}} g(u) u^{-\alpha/2} \int_{\mathbb{R}} \exp\left\{ -A\left( \sigma_+(\hat{s}) + \frac{1}{2}\sigma_+''(\hat{s})s^2/u^\alpha \right)^\alpha u^\alpha \right\} ds$$
$$\sim \frac{\sqrt{2\pi} H\sigma_+(\hat{s})^\rho}{\hat{s}\sqrt{A\alpha\sigma_+''(\hat{s})/\sigma_+(\hat{s})^{1-\alpha}}} g(u) u^{-\alpha/2} e^{-A\sigma_+(\hat{s})^\alpha u^\alpha} \qquad \text{as } u \to \infty.$$

(4.7)

Hence, part (ii) of Corollaries 1 and 2 applies, for $\beta < 1 - \alpha$, to give the Piterbarg property for $\gamma < H + (H \wedge 1)/\alpha$, and the strong Piterbarg property for $t \geq 0$ with

$$\lim_{u \to \infty} \frac{t(u)}{u^{1/\gamma - \mathfrak{a}(1 - H/\gamma)/(H \wedge 1)}} = 0 \qquad \text{for some } \mathfrak{a} > \alpha.$$

EXAMPLE 5 —Part II.   Here we continue the study of the case when $r$ is given by (4.6), in the particular case when $f = \mathbf{1}_{(0,1]}$ (so that $X$ has independent incre- ments). We show that, in this case, the Piterbarg property is absent, if $\gamma \geq H + 1/\alpha$.

By Theorem 1 and (3.24), the Piterbarg property is absent when

$$(4.8) \qquad \limsup_{u \to \infty} \frac{1}{R(u)} \mathbf{P}\left\{ \sup_{0 \leq r \leq t/u^{1/(\gamma - H)}} \sup_{s \geq r} \frac{X(s) - X(r)}{1 + c(s - r)^{\gamma}} > u \right\} > 1.$$

Theorem A does not apply here, since suprema are taken over regions that depend on $u$. However, the arguments for that theorem in Rosiński and Samorodnitsky (1993) produce an asymptotic lower bound, for the probability in (4.8), which implies the following sufficient condition for (4.8):

$$\limsup_{u \to \infty} \frac{1}{R(u)}$$
$$\times \nu\left\{ y \in \mathbb{R}^{(0,\infty) \cap \mathbb{Q}} : \sup_{r \in (0, t/u^{1/(\gamma - H)}) \cap \mathbb{Q}} \sup_{s \in (r, \infty) \cap \mathbb{Q}} \frac{y(s) - y(r)}{1 + c(s - r)^{\gamma}} > u \right\}$$
$$> 1.$$
$$(4.9)$$

Denoting the numerator in (4.9) by $R_t(u)$, we have, by the inequalities $\gamma \geq H + 1/\alpha$ and $(1 - x)^H \leq 1 - (H \wedge 1)x$ for $x \in [0, 1]$, together with (4.7) [cf. (4.3)],

$$R_t(u) = H \int_0^\infty \mu\left( \left\{ x > 0 : y^H x \sup_{0 \leq r \leq t/u^{1/(\gamma - H)}} \sup_{s \geq r} \frac{\mathbf{1}_{(0,s]}(y) - \mathbf{1}_{(0,r]}(y)}{1 + c(s - t)^{\gamma}} > u \right\} \right) \frac{dy}{y}$$

$$= H \int_0^{t/u^{1/(\gamma - H)}} r\left( y^{-H} u \right) \frac{dy}{y}$$

$$+ H \int_{t/u^{1/(\gamma - H)}}^\infty r\left( \frac{1 + c(y - t/u^{1/(\gamma - H)})^\gamma}{y^H} u \right) \frac{dy}{y}$$

$$\geq H \int_{t/u^{1/(\gamma - H)}}^\infty r\left( \sigma(y - t/u^{1/(\gamma - H)}) \left( \frac{y - t/u^{1/(\gamma - H)}}{y} \right)^H u \right) \frac{dy}{y}$$

$$\geq H \int_{t/u^{1/(\gamma - H)}}^\infty r\left( \sigma(y - t/u^{1/(\gamma - H)}) \left( 1 - (H \wedge 1) \frac{t/u^\alpha}{y} \right) u \right) \frac{dy}{y}$$

$$\sim \frac{H\sigma(\hat{s})^{\rho}}{\hat{s}} g(u) u^{-\alpha/2}$$

$$\times \int_{\mathbb{R}} \exp\left\{-A\left(\sigma(\hat{s}) + \frac{1}{2}\sigma''(\hat{s})s^2/u^{\alpha}\right)^{\alpha}\left(u^{\alpha} - (H \wedge 1)\frac{t}{\hat{s}}\right)\right\} ds$$

$$\sim R(u)\exp\left\{A\sigma(\hat{s})^{\alpha}(H \wedge 1)\frac{t}{\hat{s}}\right\} \qquad \text{as } u \to \infty.$$

This gives (4.9). [If nervous about this calculation, shrink the domain of integration from $[t/u^{1/(\gamma-H)}, \infty)$ to $\hat{s} \pm K/u^{\alpha/2}$, and send $K \uparrow \infty$ at the end.]

### 4.2. $\alpha$-stable processes.

First we consider a storage process $Y$, with an $\alpha$-stable $H$-s.s. input process $X$.

EXAMPLE 6. Let $X$ be a strictly $\alpha$-stable $H$-s.s. process, that satisfies Condition $\mathfrak{X}$, and is given by (2.7), where $M$ is a strictly $\alpha$-stable random measure ($S\alpha S$ if $\alpha = 1$). By calculations similar to those in Examples 3 and 4, we have

$$R(u) = 1 \wedge u^{-\alpha} \int_S \left(\frac{1-\beta}{2} \sup_{t \in \mathbb{Q}^+} \frac{(f_t(s)^-)^{\alpha}}{(1+ct^{\gamma})^{\alpha}}\right.$$

$$\left. + \frac{1+\beta}{2} \sup_{t \in \mathbb{Q}^+} \frac{(f_t(s)^+)^{\alpha}}{(1+ct^{\gamma})^{\alpha}}\right) w_0(s)\, d\lambda(s).$$

We assume that the above integral is nonzero, so that $R$ is not identically zero.

Provided that $X$ satisfies (3.23), Corollary 1 now gives the Piterbarg property for any $\gamma > H$, for the storage process $Y$, while Corollary 2 gives the strong Piterbarg property for $t(u) = o(u^{1/\gamma})$. However, by Samorodnitsky (1988) [see also Samorodnitsky and Taqqu (1994), Theorem 10.5.1], the limit in (3.23) is

$$\left(\int_S \left(\frac{1+\beta}{2} \sup_{t \in (0,1) \cap \mathbb{Q}} (f_t^-)^{\alpha} + \frac{1-\beta}{2} \sup_{t \in (0,1) \cap \mathbb{Q}} (f_t^+)^{\alpha}\right) w_0\, d\lambda\right)$$

$$\times \left(\int_S \left(\frac{1-\beta}{2} \sup_{t \in (0,1) \cap \mathbb{Q}} (f_t^-)^{\alpha} + \frac{1+\beta}{2} \sup_{t \in (0,1) \cap \mathbb{Q}} (f_t^+)^{\alpha}\right) w_0\, d\lambda\right)^{-1}.$$

This ratio is finite, by the local boundedness of $X$ and the assumption that $R$ is nonzero.

Next we consider the Piterbarg properties (1.2), in the case when the process $Y$ itself is $\alpha$-stable. Now $Y$ is no longer a storage process, and the example simply is to illustrate how unusual the Piterbarg property is for "usual processes."

EXAMPLE 7. For $Y = \{Y(s)\}_{s \geq 0}$ an $\alpha$-stable process, $\alpha \in (0, 2]$, we may write $Y = Y_1 - Y_2 + \mu$, where $\mu \colon [0, \infty) \to \mathbb{R}$ is a suitable function, while $\{Y_1(s)\}_{s \geq 0}$ and $\{Y_2(s)\}_{s \geq 0}$ are independent $\alpha$-stable processes, such that (for $j = 1, 2$ and $s \geq 0$)

$$\mathbf{E}\{e^{i\theta Y_j(s)}\} = \begin{cases} \exp\left\{-\sigma^{\alpha}_{Y_j(s)}|\theta|^{\alpha}\left(1 - i\tan\left(\frac{\pi\alpha}{2}\right)\mathrm{sign}(\theta)\right)\right\}, & \text{if } \alpha \neq 1, \\ \exp\left\{-\sigma_{Y_j(s)}|\theta|\left(1 + i\frac{2}{\pi}\mathrm{sign}(\theta)\ln(|\theta|)\right)\right\}, & \text{if } \alpha = 1, \end{cases}$$

for $\theta \in \mathbb{R}$. Here $\sigma_{Y_j(s)}$ is the scale parameter of the $\alpha$-stable random variable $Y_j(s)$. In the Gaussian case $\alpha = 2$, we may take $Y_1 = 0$. We assume that $Y$ satisfies Condition $\mathfrak{X}$, from which it follows that $Y_1$, $Y_2$ and $\mu$ can be taken to satisfy Condition $\mathfrak{X}$.

We are going to investigate when the Piterbarg property (1.2a) holds.

CASE 1. If $Y_1 \overset{d}{\neq} 0$, then (1.2a) holds if and only if the f.d.d.'s of $Y_1$ coincide with those of a single $\alpha$-stable random variable [Samorodnitsky and Taqqu (1994), Theorem 10.5.1].

CASE 2. If $Y_1 \overset{d}{=} 0$ and $\alpha \geq 1$, then (1.2a) holds if and only if the f.d.d.'s of $Y_2$ coincide with those of a single $\alpha$-stable random variable, and $\mu$ is constant.

To see this, notice that $\sigma_{Y_2}$ and $\mu$ must be constants on $[0, t]$ [Samorodnitsky and Taqqu (1994), equation 1.2.11]. Given these properties, we have

$$(4.10) \quad \begin{aligned} \mathbf{P}&\left\{\sup_{s \in [0, t]} Y(s) > u\right\} \\ &\geq \mathbf{P}\{\{Y(r) > u\} \cup \{Y(s) > u\}\} \\ &\geq 2\mathbf{P}\{Y(r) > u\} - \mathbf{P}\{\tfrac{1}{2}(Y(r) + Y(s)) > u\} \end{aligned}$$

for $r, s \in [0, t]$. If $\alpha > 1$, then the second probability on the right-hand side is $o(\mathbf{P}\{Y(r) > u\})$, unless $\sigma_{Y(r)+Y(s)} = \sigma_{Y(r)} + \sigma_{Y(s)}$. By Minkowski's inequality, this happens if and only if $Y(r) = Y(s)$ a.s. If $\alpha = 1$, then the spectral measure $\Gamma$ of $(Y(r), Y(s))$ [Samorodnitsky and Taqqu (1994), Section 2.3] is supported on $S_2^- = \{(s_1, s_2) \in \mathbb{R}^2 : s_1^2 + s_2^2 = 1, s_1, s_2 \leq 0\}$, and [Samorodnitsky and Taqqu (1994), Example 2.3.4]

$$\begin{aligned} &\frac{Y(r) + Y(s)}{2} \\ &\overset{d}{=} Y(r) - \frac{2}{\pi}\left(\int_{S_2^-} \frac{s_1 + s_2}{2}\ln\left|\frac{s_1 + s_2}{2}\right| d\Gamma(s) - \int_{S_2^-} s_1 \ln|s_1| \, d\Gamma(s)\right). \end{aligned}$$

By a convexity argument, unless $Y(r) = Y(s)$ a.s., the term $\frac{2}{\pi}(\cdots)$ on the right-hand side is strictly positive, so that the second probability on the right-hand side in (4.10) is $o(\mathbf{P}\{Y(r) > u\})$ [Samorodnitsky and Taqqu (1994), equation 1.2.12].

See Talagrand (1988) and Albin (1999) for more information related to Case 2.

CASE 3. If $Y_1 \overset{d}{=} 0$ and $\alpha < 1$, then (1.2a) holds (in the sense of $0/0 = 1$) since $Y_2$ is nonnegative.

SPECIAL CASE. If $Y$ is $S\alpha S$, then (1.2a) holds if and only if the f.d.d.'s of $Y$ coincide with those of a single $\alpha$-stable random variable.

Turning to (1.2b), with $t = t(u) \to \infty$ as $u \to \infty$, the above characterizations remain valid [with (1.2b) replacing (1.2a)], if appropriate global boundedness properties are imposed on $Y_2$ and $\mu$ in Case 1, and on $\mu$ in Case 3.

**Acknowledgment.** We are grateful to an anonymous referee for very useful comments.

DEPARTMENT OF MATHEMATICS
CHALMERS UNIVERSITY OF TECHNOLOGY
412 96 GÖTEBORG
SWEDEN
E-MAIL: palbin@math.chalmers.se

SCHOOL OF OPERATIONS RESEARCH
  AND INDUSTRIAL ENGINEERING
  AND DEPARTMENT OF STATISTICAL
  SCIENCE
CORNELL UNIVERSITY
ITHACA, NEW YORK 14853
USA
E-MAIL: gennady@orie.cornell.edu